\documentclass[11pt]{amsart}

\usepackage{amssymb}
\usepackage{graphics}
\usepackage{graphicx}
\usepackage{amscd}

\newtheorem{theorem}{Theorem}
\newtheorem{lemma}{Lemma}

\theoremstyle{definition}
\newtheorem{remark}{Remark}

\newcommand{\HF}{\widehat{HF}}
\newcommand{\on}{\operatorname}

\renewcommand{\d}{\partial}

\newcommand{\h}{\widehat}

\newcommand{\x}{{\bf x}}
\newcommand{\y}{{\bf y}}
\newcommand{\z}{{\bf z}}
\newcommand{\co}{{\bf c}}

\newcommand{\Z}{\mathbb{Z}}

\newcommand{\D}{\mathbb{D}}
\newcommand{\T}{\mathbb{T}}

\newcommand{\Sym}{\on{Sym}}
\renewcommand{\Re}{\on{Re}}

\begin{document}

\author{Olga Plamenevskaya}
\address{Department of Mathematics, SUNY Stony Brook, Stony Brook, NY 11794}
\email{olga@math.sunysb.edu}
\title{A combinatorial description \\ of  the Heegaard Floer contact invariant}

\begin{abstract} In this short note, we observe that the Heegaard Floer 
contact invariant is combinatorial by applying  the algorithm  of Sarkar--Wang 
to the description of the contact invariant due to Honda--Kazez--Mati\'c.
We include an example of this combinatorial calculation. 
\end{abstract}

\maketitle

\section{Introduction}

Recent months have seen a significant advance in Heegaard Floer theory: it turned 
out  that certain Heegaard Floer 
homologies admit a purely combinatorial description. In particular, it was shown in  
\cite{MOS} that the Heegaard Floer homologies of a knot can 
be computed from a grid diagram of the knot by a simple combinatorial procedure.
 Heegaard Floer homology $\HF(Y)$
of a 3-manifold  $Y$ also admits a combinatorial description \cite{SW}, but this description is 
less straighforward. Starting from an arbitary admissible Heegaard diagram of $Y$,
one has to change   the $\beta$ curves by isotopies and handleslides so that 
in the resulting Heegaard diagram almost all the domains of holomorphic disks are squares 
or bigons. It is then easy to understand the moduli spaces of the holomorphic disks 
needed to compute the differential, since squares and bigons with Maslov index =1 admit a unique   
holomorphic representative. (We assume that the reader is familiar with the basic setup
of the Heegaard Floer theory; see \cite{OSs} for a survey.)

For a contact 3-manifold $(Y, \xi)$,  %with $\xi$ a co-orientable contact structure.
Ozsv\'ath and Szab\'o introduce an invariant $c(\xi)$   which is a distinguished 
element of $\HF(-Y)$  (defined up to sign for the theory with $\Z$ coefficients)\cite{OS}.
  Since non-vanishing of $c(\xi)$ implies that $\xi$ is tight, 
this invariant gives a powerful tool  for establishing the 
tightness of a contact structure \cite{LS}.  
  
The invariant $c(\xi)$ is defined in \cite{OS} via an open 
book decomposition of $(Y, \xi)$, as follows. If the genus of the page of the open book is $g$ 
(and the binding is connected), we consider a certain  Heegaard diagram 
for  $Y$ of genus $2g+1$, compatible with the open book.  
The generators of the group $\h{CF}(-Y)$ are given 
by $(2g+1)$-tuples of  the intersection points of $\alpha$- and $\beta$-curves in the diagram,
and a distinguished $(2g+1)$-tuple gives a cycle  $\co$ which descends to the invariant $c(\xi)$ 
in homology.    
    
\begin{figure}[ht]
\includegraphics[scale=0.7]{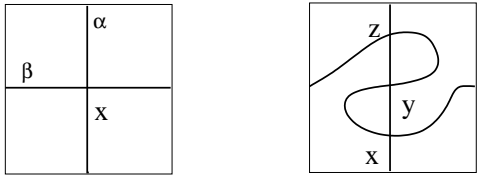}
\caption{The Heegaard surface here is a torus obtained by identifying the opposite sides of the square. After the isotopy of the $\beta$-curve,  $\x$ is no longer a cycle in the Heegaard Floer complex.}
\label{isotopies}
\end{figure}

Our goal is to show that $c(\xi)$ can be computed in a combinatorial fashion.  We would like to apply Sarkar--Wang algorithm \cite{SW} to obtain a Heegaard diagram where the holomorphic disks can be easily identified. However, 
we are concerned with a specific cycle, not the homology group as a whole, and the 
isotopies and/or handleslides of the $\beta$-curves  performed  
on the Heegaard surface would possibly affect the generator $\co$.
Indeed, the homology class of a geometric generator $\co$ can change even if all the isotopies 
are supported away from the intersection points forming $\co$;  $\co$ might even  no longer 
be a cycle after the isotopy. 
Consider  for example the 
genus 1 Heegaard diagram for $S^1$ given by an $\alpha$ and a $\beta$ curve on the torus intersecting 
at one point. This intersection point $\x$ is a cycle which  generates  $\HF(S^3)=\Z$. 
We isotope
the $\beta$-curve to introduce two extra intersection points, $\y$ and $\z$, as shown on Figure \ref{isotopies},
Now, we have $\d \x=\y$, $\d \z=\y$ in the chain complex $\h{CF}(S^3)$ for the new Heegaard diagram;
 so $\x$ is no longer a cycle, and    $\HF(S^3)$ is generated by $\x-\z$.

This shows that we need to be more careful.   Fortunately,  there is an alternate 
geometric description of $c(\xi)$ due to Honda--Kazez--Mati\'c: in \cite{HKM}, 
they show that $c(\xi)$ can be found from a Heegaard diagram which is different from and somewhat 
simpler that the one  in \cite{OS}.
This alternative Heegaard diagram is  easier to handle, and fits well with the Sarkar--Wang algorithm.

We combine the arguments from \cite{HKM} and \cite{SW} to establish our main result in the next section.

\section{The Main Result}

We first recall the construction from \cite{HKM}. Let $(S, h)$ be an open book 
decomposition for the contact manifold $(Y,  \xi)$; here $S$ denotes the page 
of the open book, and $h$ the monodromy. This means that $Y$ is homeomorphic 
to $S\times [0,1]/\sim$, where the equivalence relation $\sim$ is given by
\begin{align*} (x, 1) &\sim (h(x),0), \qquad  x\in S \\
               (x, t) & \sim (x, t'), \qquad  x\in \d S, t, t'\in [0,1].
\end{align*}
 The open book produces a Heegaard splitting  $Y=H_1\cup H_2$, with 
 $H_1=  S\times [0,1/2]/\sim$, $H_2=S\times [1/2,1]/\sim$. The Heegaard diagram for $Y$
 can then be given by the Heegaard surface $\Sigma = S_{1/2}\cup -S_0$, and
 the $\alpha$- and $\beta$-curves, defined as follows. Consider a set of disjoint, properly embedded
 arcs $\{a_1, a_2, \dots, a_n\}$ on $S$ such that $S\setminus \bigcup a_i$ is a single polygon. Obtain arcs $b_i$ by 
 changing the arcs $a_i$ via a small isotopy so that the endpoints of $a_i$ are isotopied along $\d S$
 (in the direction dictated by the boundary orientation), the arcs $a_i$ and $b_i$ intersect transversely 
 at one point, and the sign of this intersection in positive (the orientation of $b_i$ is induced 
 from the orientation of $a_i$ by the isotopy). The curves $\alpha_i= \d(a_i \times [0, 1/2])$ and 
 $\beta_i = \d(b_i \times [1/2, 1])$ form attaching circles for the handlebodies $H_1$ and $H_2$; 
 they can be thought of as the $\alpha$- and $\beta$- curves on $\Sigma$. We can write 
  $$
  \alpha_i= a_i\times \{ 1/2\} \cup a_i\times \{0\}, \qquad \beta_i= b_i\times \{ 1/2\} \cup h(b_i)\times \{0\};
  $$
 thus, the intersection of  $\alpha$- and $\beta$- curves with $S_{1/2}$ is completely standard 
 (and given by $a_i$, $b_i$); the picture on  $S_{0}$  depends on the monodromy $h$. 
 For an illustration of such Heegaard diagram, see \cite{HKM}, or look at the example in the next section.
 The basepoint $z_0$
 is placed  on $S_{1/2}$ in the polygonal region (not in the thin strips between $a_i$'s and $b_i$'s); we denote 
 this polygonal region by $D_0$.
 Now, let $c_i$ be the intersection point between $a_i$ and $b_i$ on $S_{1/2}$. It is shown in \cite{HKM}
 that the element $\co =(c_1, c_2, \dots, c_n)\in \h{CF}(\Sigma, \beta, \alpha, z_0)$ is a cycle which 
 descends to the element $c(\xi)$ in homology $H_*(\h{CF}(\Sigma, \beta, \alpha, z_0))=\HF(-Y)$. 
     (Note that the roles of 
the $\alpha$- and $\beta$-curves are interchanged because we need the homology of $-Y$ instead 
of the homology of $Y$. So for a Heegaard diagram of genus $g$, 
a Whitney disk from $\x$ to $\y$ is now a map $\phi:\D\to \Sym^g \Sigma$ such that $\phi(i)=\x$, $\phi(-i)=\y$, and $\D\cap \{\Re z>0\}$ is mapped into $\T_\alpha=\alpha_1\times \dots \times \alpha_g$, while
 $\D\cap\{ \Re z<0\}$ is mapped into $\T_\beta=\beta_1 \times \dots \times \beta_g$. 
 This does not affect the combinatorial algorithm that we will be using later.)

Applying the idea of Sarkar--Wang, we would like to find an open book decomposition for $Y$ 
such that all but one regions in the corresponding Heegaard diagram are bigons or squares (here and below, a region is 
a connected component in the complement of $\alpha$- and $\beta$-curves in $\Sigma$). 

\begin{theorem} \label{mythm} There exists an open book $(S, h')$ for $(Y, \xi)$, such that the Heegaard diagram 
described above has only disk and square regions (except for the polygonal region $D_0\subset S_{1/2}$). The   
monodromy $h'$ differs from the monodromy of given open book $(S, h)$ by an isotopy, i.e. 
$h'=\phi \circ h$, where $\phi:S \to S$ is a diffeomorphism fixing the boundary and isotopic to identity. 
\end{theorem} 

\begin{proof} The algorithm of \cite{SW} tells us to get rid of non-disk regions and $2n$-gons 
with $n>2$ by performing  isotopies (finger moves) of $\beta$-curves. First, the non-disks are killed;
then, the $2n$-gons are dealt with (one after another) roughly as follows. We look at the ``distance'' between a given region $D$ and the region $D_0$ (the minimal number of intersections between the $\beta$-curves and an arc connecting 
$z_0$ to an interior point of $D$), and number all the regions $D_0$, $D_1$, $D_2$, etc, so that 
the distance between $D_m$ and $D_0$ increases with $m$. Sarkar and Wang explain how to perform finger moves 
that break up  a given  $2n$-gon $D_k$ into polygons with fewer sides, pushing part of the boundary of $2n$-gon into 
other regions (typically labelled with  smaller numbers). During  this process, the regions $D_l$ 
with $l>k$ which are already bigons or squares remain bigons or squares, so the process eventually terminates.

We observe that for our Heegaard diagram coming from an open book, all the finger moves can be performed 
in the $S_0$ part of $\Sigma$, since every boundary curve of every region in $\Sigma$ intersects $S_0$.

Apart from the finger moves, handleslides of the $\beta$-curves are sometimes needed in the algorithm of \cite{SW}. 
After such a handleslide, the Heegaard diagram would no longer be compatible with considerations in \cite{HKM}.
However, it turns out that the need for handleslides fortunately does not arise in our case. 
Indeed, a handleslide in \cite{SW} is only needed when a finger is pushed through a collection of adjacent regions, 
none of which has a smaller distance from $D_0$ (in particular, none of these regions is $D_0$), 
and then comes back to the region where it started. This means that the finger goes around a full copy 
of some curve $\beta_i$. Because each of the $\beta$-curves forms part of the boundary of $D_0$,  
and $a_i$ and $b_i$ intersect at one point in $S_{1/2}$, it follows that the finger has to go through $D_0$, 
which is a contradiction. 

Therefore, we can obtain  a ``nice'' Heegaard diagram by performing a sequence of isotopies on $S_0\subset \Sigma$
(away from boundary). A composition of these isotopies gives a diffeomorphism $\phi$ and the open book 
$(S, \phi\circ h)$ equivalent to the open book $(S, h)$ we started with. \end{proof}

As shown in \cite{SW}, finding the differentials in the ``hat'' Heegaard Floer complex for the diagram that we 
obtained is now a combinatorial matter. Let $g$ be the genus of our Heegaard diagram; 
stabilizing the open book if necessary, we can assume that $g\geq 2$. We consider a generic almost complex structure on $\Sym^g \Sigma$
which is a small perturbation of the product complex structure,
 and look for holomorphic representatives 
of a Whitney disk  $\phi \in \pi_2(\x, \y)$ connecting 
the points $\x$ and $\y$; the domain of $\phi$ is then a linear combination of some regions $D_i$,
$D=\sum_{i=1}^k  a_i D_i$. If $\phi$ has holomorphic representatives, we must have $a_i\geq 0$. 
We assume that the Maslov index of $\phi$ is 1. It remains to cite

\begin{lemma}(\cite{SW}) Under these conditions, $\phi$ has a unique holomorphic representative. Moreover, 
$D$ is either an embedded square tiled by squares, or a bigon tiled by squares and a bigon.
\end{lemma}

We conclude this section with the following 

\begin{remark}
If we only want to find out whether the invariant $c(\xi)$  vanishes or not (without pinpointing 
it in the homology group $\HF(-Y)$),  it often suffices to simplify the regions that contribute to the moduli 
spaces for the differentials possibly killing  the element $\co$. These are the regions passing through 
the ``thin strips'' on $S_{1/2}\subset \Sigma$; all of them are adjacent to $D_0$, so they can be converted 
into bigons or squares simply by pushing parts of their boundary into $D_0$. If we find any points
$\x$ such that the expression for $d\x$ contains $\co$, we will need to understand the full boundary of $\x$ (and so to convert more regions into bigons or squares); however, all the other 
 regions in the Heegaard diagram may be left as is.       
\end{remark}

\section{An Example}
\begin{figure}[ht]
\includegraphics[scale=0.7]{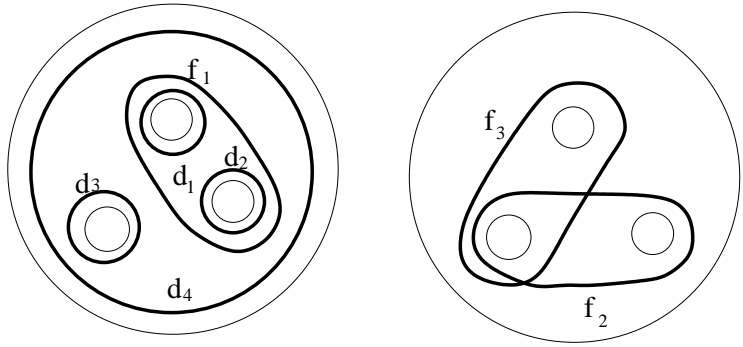}
\caption{The monodromy of the open book is given by the product of Dehn twists around the pictured curves.}
\label{lantern}
\end{figure}

In this section we illustrate the combinatorial calculation of 
$c(\xi)$ by
 the following example. (We consider $\Z/2$ coefficents for simplicity, and only 
 show that $c(\xi)$  is non-zero without computing $\h{HF}(-Y)$). Consider the contact manifold $(Y, \xi)$
given by  the open book whose page is a four-punctured 
sphere, and the monodromy is  the composition of  positive Dehn twists around the boundary 
curves $d_1$, $d_2$, $d_3$, $d_4$ and a negative Dehn twist around the curve $f_1$ (see Figure \ref{lantern}). (Since 
the curves are disjoint, the order in which the Dehn twists are performed is not important.)
This open book is easy to understand without Heegaard Floer theory: using the lantern relation \cite{D}, we can express the monodromy 
as the product $T_{f_2}^+ T_{f_3}^+$ of the positive Dehn twists around the curves $f_2$ and $f_3$.
By \cite{Gi}, this shows that $\xi$ is Stein fillable (and therefore tight); moreover, we can perform 
two positive destabilizations to see  that $\xi$ is in fact simply the standard tight contact 
structure on $S^1\times S^2$ (given by the open book with an annular page and trivial monodromy).

\begin{figure}[ht]
\includegraphics[scale=0.7]{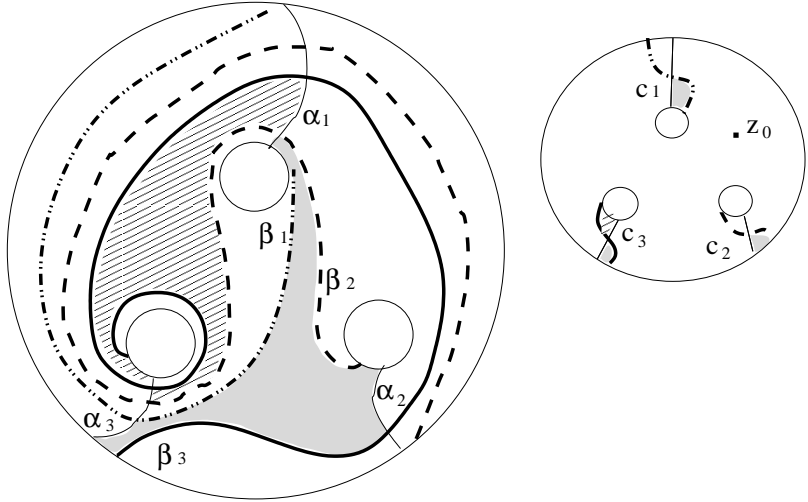}
\caption{This is a Heegaard diagram for $(Y, \xi)$.  The surfaces $S_0$ and $S_{1/2}$, the "bottom" and "top" parts of $\Sigma$, are shown separately    ($S_{1/2}$ is smaller because it's completely 
standard).  Note that $\Sigma= -S_0\cup S_{1/2}$; the picture shows $S_0$, not $-S_0$. The thin curves are the $\alpha$-curves, the thicker lines (solid, dashed and dash-dotted)
are the $\beta$-curves.  The contact element $\co=(c_1, c_2, c_3)$ lies on  $S_{1/2}$. The ``bad'' regions on $\Sigma$ are shown.}
\label{shadediag}
\end{figure}

\begin{figure}[ht]
\includegraphics[scale=0.7]{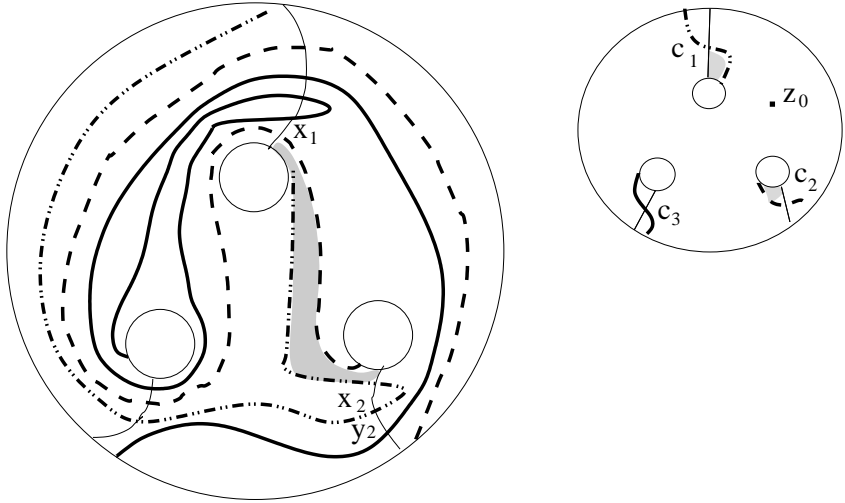}
\caption{The Heegaard diagram after winding. The domain of the holomorphic disk connecting 
$\x=(x_1, x_2, c_3)$ and $\co=(c_1, c_2, c_3)$ is shaded.}
\label{winddiag}
\end{figure}

Our point, however, is to apply Theorem \ref{mythm} and to perform a calculation for  the original open book. 
First, we look at the Honda--Kazez--Mati\'c-style
 Heegaard diagram for $(Y, \xi)$ (Figure \ref{shadediag}).  Observe that there are two ``bad'' regions: a non-disk region
 (hatched)  and a hexagonal 
 region (shaded) in the complement of $\alpha$ and $\beta$ curves 
 on the Heegaard surface. 
 We get rid of them by applying the Sarkar--Wang algorithm as described in Theorem \ref{mythm}, which in this case amounts simply 
 to winding two of  the $\beta$-curves as shown in Figure \ref{winddiag}: as explained in the end 
 of section 2, we just have to push fingers out of bad regions into $D_0$.  
 
 We obtain a Heegaard diagram shown on Figure \ref{winddiag}.  Examining it, we see that there is 
 only one possible domain of a differential going to $\co$ from another point $\x$; the point $\x=(x_1, x_2, c_3)$
 is shown in the picture, and the domain is shaded. (Recall that the intersection of all such domains with $S_{1/2}$ should lie in 
 the thin strips between $a_i$ and $b_i$, since the domain must not contain $z_0$.) But then we have 
 $d \x= \co+\y$, where $\y=(x_1, y_2, c_3)$, since there is a bigon connecting $\x$ and $\y$
 (and no other Whitney disks from $\x$). This shows 
 that $\co$ is not a boundary, and so $c(\xi)\neq 0$. We can conclude that the contact structure $\xi$
 is tight.    


\begin{thebibliography}{WWW}
 
 \bibitem{HKM} K. Honda, W. Kazez, G. Mati\'c, {\em On the contact class in Heegaard Floer homology},  math.GT/0609734.
 
 \bibitem{D} M. Dehn, {\em Papers on group theory and topology}, Springer--Verlag, New York, 1987.
 
 \bibitem{Gi} E.~Giroux, 
{\em G\'eom\'etrie de contact: de la dimension trois vers les 
dimensions sup\'erieures}, in: {\em Proceedings of the International Congress of Mathematicians, Vol. II (Beijing, 2002)},  Higher Ed. Press, Beijing, 2002, 405--414.
 
 \bibitem{LS}  P.~Lisca and A.~Stipsicz, {\em Heegaard Floer 
Invariants and Tight
 Contact Three--Manifolds I-III}, Geom. Topol. {\bf 8} (2004) 925-945,  math.SG/0404136,
 math.SG/0505493.
 
 \bibitem{MOS} C. Manolescu, S. Sarkar, P. Ozsv\'ath, {\em  On combinatorial link Floer homology}, 
math.GT/0610559. 
 
 \bibitem{OSs} P.~Ozsv\'ath and Z.~Szab\'o, {\em Heegaard diagrams and Floer homology}, math.GT/0602232.
 
 \bibitem{OS} P.~Ozsv\'ath and Z.~Szab\'o, {\em Heegaard Floer
  homologies and contact structures},  Duke Math. J. {\bf 129}   (2005),  no. 1, 39--61.
 
 \bibitem{SW} S. Sarkar, J. Wang,  {\em A combinatorial description of some Heegaard Floer homologies}, math.GT/0607777.
 
 
\end{thebibliography}
\end{document}